\documentclass[12pt]{article}
\usepackage{amsfonts,mathrsfs,amssymb}
\def\d{{\rm d}}
\def\mi{{\rm i}}
\def\eps{\varepsilon}
\def\Re{\mathop{\rm Re\,}\nolimits}
\def\Im{\mathop{\rm Im\,}\nolimits}
\def\e{\mathop{\rm e}\nolimits}
\def\Res{\mathop{\rm Res}\nolimits}
\def\hf{{\textstyle{1 \over 2}}}
\def\qt{{\textstyle{1 \over 4}}}

\title{Sharpenings of Li's criterion\\
for the Riemann Hypothesis}

\author{{\bf Andr\'e Voros}\footnote{Also at: 
Institut de Math\'ematiques de Jussieu--Chevaleret (CNRS UMR 7586), 
Universit\'e Paris 7, F-75251 Paris CEDEX 05, France.}\\
\\
CEA, Service de Physique Th\'eorique de Saclay\\
(CNRS URA 2306)\\
F-91191 Gif-sur-Yvette CEDEX (France)\\
E-mail : {\tt voros@spht.saclay.cea.fr}}

\begin{document}

\maketitle

{\abstract
Exact and asymptotic formulae are displayed for the coefficients
$\lambda _n$ used in Li's criterion for the Riemann Hypothesis.
For $n \to \infty$ we obtain that if (and only if) 
the Hypothesis is true, $\lambda _n \sim n(A \log n +B)$
(with $A>0$ and $B$ explicitly given, 
also for the case of more general zeta or $L$-functions);
whereas in the opposite case, $\lambda _n$ has a non-tempered oscillatory form.
}
\bigskip

\emph{Li's criterion} for the Riemann Hypothesis (RH) states that the latter 
is true if and only if a specific real sequence $\{ \lambda_n \}_{n=1,2,\ldots}$
has \emph{all its terms positive} \cite{LI1,BL}.
Here we show that it actually suffices to probe the $\lambda _n$ for their 
\emph{large-$n$ behavior}, which fully encodes the Riemann Hypothesis 
by way of a clear-cut and explicit \emph{asymptotic alternative}.
To wit, we first represent $\lambda _n$ exactly by a finite oscillatory sum (\ref{LZS}),
then by a derived integral formula (\ref{LZI}), which can finally be
evaluated by the saddle-point method in the $n \to +\infty$ limit.
As a result, $\lambda _n$ takes one of two
sharply distinct and mutually exclusive asymptotic forms:
if RH is true, $\lambda _n$ will \emph{grow tamely}
according to~(\ref{RER}); if RH is false,
$\lambda _n$ will oscillate with an \emph{exponentially growing amplitude},
in both $+$ and $-$ directions, as described by~(\ref{RS}).
This dichotomy thus provides a sharp criterion of a new \emph{asymptotic} type 
for the Riemann Hypothesis 
(and for other zeta-type functions as well, replacing (\ref{RER}) by (\ref{REs})).

This work basically reexposes our results of April 2004 announced in \cite{V}, 
but with an uncompressed text; 
we also update the references and related comments: for instance, we now derive
as (\ref{TS}) a large-$n$ expansion surmised by Ma\'slanka \cite{M1} in the meantime.

\section{Background and notations.}

We study the sequence \cite{K,LI1}
(in the notations of Li, whose $\lambda_n$ are $n$ times Keiper's)
\begin{equation}
\label{LDef}
\lambda_n = \sum_\rho \, [1-(1-1/\rho)^n] \qquad (n=1,2, \ldots),
\end{equation}
where $\rho$ are the nontrivial zeros of Riemann's $\zeta(s)$,
grouped by pairs within summations and products as
\begin{equation}
\label{ZER}
\{\rho = \hf \pm \mi \tau_k \}_{k=1,2,\ldots}, \qquad
\Re \tau_k \mbox{ positive and non-decreasing;}
\end{equation}
we also parametrize each such pair by the single number
$x_k=\rho(1-\rho)=\qt+{\tau_k}^2$.

We will use the completed zeta function $\Xi(s)$ (normalized as $\Xi(0)=\Xi(1)=1$)
and a \emph{symmetrized} form of its Hadamard product formula \cite{Ed,Vz},
\begin{equation}
\label{XI}
\Xi(s) = s(s-1) \Gamma (s/2) \pi^{-s/2} \zeta(s) 
\equiv \prod_{k=1}^{\infty} \, \Biggl[ 1 - {s(1-s) \over x_k} \Biggr] ;
\end{equation}
we will also use  a ``secondary" zeta function built over the Riemann zeros,
\begin{equation}
\label{ZDef}
Z(\sigma) = \sum_{k=1}^\infty {x_k}^{-\sigma}, \qquad \Re \sigma > \hf ,
\end{equation}
which extends to a meromorphic function in $\mathbb C$ having all its poles
at the negative half-integers, plus one pole at $\sigma=+\hf$ \cite{Ku}
of polar part \cite{Vz}
\begin{equation}
\label{DPol}
Z(\hf+\eps) = R_{-2} \,\eps^{-2} + R_{-1} \,\eps^{-1} + O(1)_{\eps \to 0} \, , 
\end{equation}
with
\begin{equation}
\label{EX}
R_{-2} = (8\pi)^{-1} , \quad R_{-1} = - (4\pi)^{-1} \log 2\pi
\qquad \mbox{in this case.}
\end{equation}

Our results \cite{V} mainly relate to those of Keiper \cite{K},
of which we only learned later (thanks to K. Ma\'slanka; they were almost never cited),
of Bombieri--Lagarias \cite{BL} on Li's criterion \cite{LI1}, 
and of Ma\'slanka \cite{M1}.
Other earlier works considering the $\lambda _n$ are \cite{BPY,C1}.

\section{New exact forms for $\lambda _n$.}

To reexpress the $\lambda _n$, we start from their generating function \cite{K,LI1}
\begin{equation}
\label{LId}
f(z) = {\d \over \d z} \log \Xi \Bigl( {1 \over 1-z}\Bigr) \equiv
\sum_{n=1}^\infty \lambda _n z^{n-1} .
\end{equation}
Now the infinite product formula in (\ref{XI}) implies
\begin{equation}
\label{HAL}
\log \Xi \Bigl( {1 \over 1-z} \Bigr) =
\sum_{k=1}^{\infty} \log \Biggl[ 1 + {z \over (1-z)^2 x_k} \Biggr] =
- \sum_{j=1}^{\infty} {(-1)^j \over j}
\Biggl[ {z \over (1-z)^2 } \Biggr] ^j Z(j);
\end{equation}
then, expanding $(1-z)^{-2j}$ by the generalized binomial formula,
reordering in powers of $z$ and substituting the output into (\ref{LId}),
we get as first result
\begin{equation}
\label{LZS}
\lambda_n = - n \sum_{j=1}^n {(-1)^j \over j}{n+j-1 \choose 2j-1} \, Z(j) .
\end{equation}

Other sums related to the $Z(k)$ are ${\mathscr Z}_j = \sum_\rho \rho^{-j}$ 
\cite{M,Le,Vz} (often denoted $\sigma_j$, but here we use $\sigma$ as variable).
It was already known that $ \lambda_n = 
\sum_{j=1}^n (-1)^{j+1} {\textstyle n \choose \textstyle j} {\mathscr Z}_j$ 
\cite[Equation~(27)]{K}, and that the ${\mathscr Z}_j$ in turn 
are complicated polynomials in the \emph{Stieltjes constants} 
$\{ \gamma_k \}_{k<j}$ \cite{M} 
(for $\lambda_n$ and $\gamma_k$ see also \cite{M2,C1,C2} and references therein).
Now the latter relations boil down to 
${\mathscr Z}_j = 1 - (1-2^{-j}) \zeta(j) +(-1)^j \eta_{j-1}$ 
\cite[Equation~(46)]{Vz}
simply by promoting \emph{logarithmic} coefficients $\eta_j$
\cite{I,BL} (cf. also \cite[Equation~(12)]{Le})
\begin{equation}
\label{SC}
\log [s \zeta(1+s)] \equiv - \sum_{n=1}^\infty \eta_{n-1} {s^n \over n} 
\end{equation}
in place of
\[
s \zeta(1+s) = 1 - \sum_{n=1}^\infty \gamma_{n-1} {(-s)^n \over (n-1)!} \, .
\]
The $\lambda_n$ thus express as affine combinations of the $\eta_j$
\cite[thm~2]{BL}. 

Remarks: 

- the $\eta_j$ are the \emph{Stieltjes} [constants'] \emph{cumulants}, 
up to some relabelings \cite{Vz,Vl};

- the $\eta_j$ admit an \emph{arithmetic} expression over the primes 
\cite[Equation~(4.1)]{BL};
see also \cite{H}, which cites \cite{VP} for the case $\eta_0 = -\gamma$;

- relation (\ref{LZS}) can be inverted also in closed form,
by the same technique as for \cite[Equation~(48)]{Vz}:
\begin{equation}
\label{SZL}
Z(j) = \sum_{n=1}^j (-1)^{n+1} {2j \choose j-n} \, \lambda_n \, .
\end{equation}

The expression (\ref{LZS}) for $\lambda_n$ has some distinctive advantages: 
it involves the functional equation $\Xi(s)=\Xi(1-s)$ through (\ref{XI});
and unlike the ${\mathscr Z}_j$, the $Z(j)$ are positive and gently varying factors:
the function $Z(\sigma)$ is regular and very smooth for real $\sigma \ge 1$.
Still, (\ref{LZS}) is an \emph{oscillatory} sum, hence difficult to control 
directly.
\medskip

Now an integral representation, equivalent to (\ref{LZS}) 
simply by residue calculus, will nevertheless prove much more flexible:
\begin{equation}
\label{LZI}
\lambda_n = {(-1)^n n \,\mi \over \pi} \oint_C I(\sigma) \,\d \sigma,
\qquad I(\sigma) =
{ \Gamma (\sigma+n) \Gamma (\sigma-n) \over \Gamma (2 \sigma+1) } \, Z(\sigma) ,
\end{equation}
where $C$ is a positive contour encircling just the subset of poles
$\sigma=+1, \ldots, +n$ of the integrand $I(\sigma)$.

\section{Asymptotic alternative for $\lambda _n, \ n \to \infty$.}

The integral formula (\ref{LZI}) readily suggests an asymptotic 
($n \to \infty$) evaluation by the classic \emph{saddle-point method} 
\cite{Er}, using $|I(\sigma)|$ as height function.
First, the integration contour $C$ is to be moved in the direction 
of decreasing $|I(\sigma)|$ as far down as possible: 
it will thus pass through some saddle-points $\sigma_\ast$ of $|I(\sigma)|$.
Then for large $n$, $I(\sigma)$ \emph{peaks} near each of these points $\sigma_\ast$, 
where it makes a contribution of the order of magnitude $|I(\sigma_\ast)|$
to the integral: thus the highest saddle-point(s) give(s) the dominant behavior. 
Consistent asymptotic approximations can also be made inside $I(\sigma)$ 
throughout: e.g. here, Stirling formulae used for $\Gamma(n+{\rm const})$.
\medskip

\emph{This approach, for an integrand not controlled in fully closed form, 
partly retains an experimental character}.
We currently advocate it for this problem as a heuristic, rather than rigorous, 
tool: it predicts the global structure of the results at once, 
and it treats all the cases readily and correctly, as other techniques confirm.
\medskip

In the present problem, for large $n$
the landscape of the function $|I(\sigma)|$ is dominantly controlled: 
by its $\Gamma $ factors, asymptotically 
$\sim \pi \, [\sin \pi \sigma \, \Gamma (2 \sigma+1)]^{-1} n^{2\sigma -1}$
for finite $\sigma$; and by the polar parts of $I(\sigma)$ near its poles.
The induced contour deformation starts as a dilation of $C$ away from 
the segment $[1,n]$ in all directions, and goes to infinity in the directions
$|\arg \sigma | < {\pi \over 2} - \delta$.
The encountered saddle-points can be of two types here (once $n$ is large enough).

1) For $\sigma$ on the segment $(\hf,1)$,
$ |I(\sigma)| \sim \pi \,
[\sin \pi \sigma \, \Gamma (2 \sigma+1)]^{-1} n^{2\sigma -1} Z(\sigma)$
always has one \emph{real} minimum
$ \sigma_{\rm r}(n) $ (tending to $\hf$ as $n \to \infty$),
which will be reached by the moving contour;
other real saddle-points lie below $\sigma=\hf$ and will not get reached here.

2) \emph{Complex} saddle-points may enter as well,
for which we may focus on the upper half-plane alone: 
the lower half-plane will give complex-conjugate (``c.c.") contributions.
As long as the moving contour stays inside a half-plane 
$\{ \Re \sigma > \hf + \eps \}$,
the integrand can be decomposed as $I=\sum_k I_k$ according to (\ref{ZDef});
then for each individual term and
within the Stirling approximation for the $\Gamma $-ratio, 
the saddle-point equation is
$0= {\d \over \d \sigma} \log |I_k(\sigma)| \sim 
{\log(\sigma^2-n^2)- 2 \log 2\sigma -\log x_k}$, 
yielding the saddle-point location
\begin{equation}
\label{CSP}
\sigma_k (n) =  n \, \mi \, / \, 2 \tau_k \, .
\end{equation}
Thus any zero \emph{on} the critical axis ($\tau_k$ real) yields
a purely imaginary $\sigma_k(n)$, not eligible:
it lies outside the domain of validity of (\ref{ZDef}),
and its contribution would be subdominant anyway.
So in the end, this paragraph excludes the real $\tau_k$.

\emph{The discussion then fundamentally splits depending on the presence
or absence of zeros off the critical axis.}

\subsection*{[RH false]}

If there is any zero $(\hf \pm \mi \tau_k)$ \emph{off} the critical axis,
we select $\arg \tau_k >0$ and assume the case of a simple zero for argument's sake.
Paragraph 2) above fully applies to each such zero: 
the complex saddle-point $\sigma_k(n)$ given by (\ref{CSP})
lies inside the domain of convergence $\{ \Re \sigma>\hf \}$
as soon as ${n > |\Im 1/\tau_k| ^{-1}}$, and for $n \to +\infty$ it gives 
an additive contribution $\sim \bigl[ (\tau_k+\mi/2 )/(\tau_k-\mi/2) \bigr] ^n$
(in the usual quadratic approximation of $\log I(\sigma)$ around $\sigma_k(n)$), 
which \emph{grows exponentially} in modulus and fluctuates in phase;
it will indeed exponentially dominate 
the contribution of the real saddle-point $\sigma_{\rm r}(n)$, to be computed later.
\medskip

This result can also be confirmed rigorously and more directly:
by a conformal mapping \cite{BL}, the function $f(z)$ in (\ref{LId})
has precisely the points $z_k={(\tau_k - \mi/2) (\tau_k + \mi/2)^{-1}}$
and $z_k^\ast$ as simple poles of residue 1 in the unit disk;
then a general Darboux theorem \cite[chap.~VII \S 2]{Di}
applies here to the poles with $|z_k|<1$,
implying that the Taylor coefficients of~$f$ (namely, the $\lambda _n$) 
indeed have the asymptotic form 
\begin{equation}
\label{DA}
\lambda _n \sim \sum_{\{ |z_k|<1 \}} z_k^{-n} + {\rm c.c.} 
\quad \pmod {o(\e^{\eps n}) \ \forall \eps >0}, \quad n \to \infty \, ,
\end{equation}
which now holds for multiple zeros as well, counted with their multiplicities.
Concretely, $\lambda_n$ then oscillates between 
\emph{exponentially growing} values \emph{of both signs}.

Infinitely many zeros $(\hf \pm \mi \tau_k)$ \emph{off the critical axis} 
are perfectly admissible here: their $z_k$ satisfy $|z_k|<1$, $z_k \to 1$, 
and the corresponding infinite sums $\sum_k z_k^{-n}$ 
still define valid $n \to +\infty$ asymptotic expansions. 
On the other hand, the general Darboux \emph{formula} (\ref{DA}) hopelessly breaks down 
if the infinitely many $z_k^{-n}$ have identical and dominant modulus,
which is precisely realized in the case [RH true], with all the $z_k$ on the unit circle!

\subsection*{[RH true]}

Here, Darboux's theorem only tells that 
$\lambda_n = o(\e^{\eps n}) \ \forall \eps >0$, 
but it fails to give any clue as to an explicit asymptotic equivalent for $\lambda_n$.
By contrast, the saddle-point treatment of the integral (\ref{LZI}) itself
remains thoroughly applicable. Simply now, all the $\tau_k$ are real,
$Z(\sigma)=O(Z(\Re \sigma) \, |\Im \sigma|^{-3/2})$ in $\{ \Re \sigma > \hf \}$,
and the contour $C$ can be freely moved towards the boundary 
$\{ \Re \sigma = \hf \}$
without meeting any of the $\sigma_k(n)$ (all of which are purely imaginary).
Hence the only dominant saddle-point is now $\sigma_{\rm r}(n) \in (\hf,1)$;
it is shaped by the double pole of $Z(\sigma)$ at $\hf$
(itself generated by the totality of Riemann zeros), so that
$ \sigma_{\rm r}(n) \sim \hf + {1 \over \log n}$.
This saddle-point is however non-isolated (it tends to the pole),
so the standard saddle-point evaluation using the quadratic approximation 
of $\log I(\sigma)$ around $ \sigma_{\rm r}(n) $ works very poorly.
Here, it is at once simpler and more accurate to
keep on deforming a portion of the contour $C$ nearest to $ \sigma=\hf$
until it fully encircles this pole (now clockwise),
and to note that the ensuing modifications to the integral are asymptotically 
smaller. Hence for [RH true], $\lambda_n$ is given (mod $o(n)$) by \cite{V}
\begin{eqnarray}
\label{REs}
\lambda_n \sim (-1)^n \, 2n \, \Res_{\sigma=1/2} I (\sigma) 
&=& 2 \pi n \, [ 2 R_{-2} (\psi(\hf+n) - 1 + \gamma ) + R_{-1} ] \nonumber \\
&=& 2 \pi n \, [ 2 R_{-2} (\log n - 1 + \gamma ) + R_{-1} ] \quad {}+ O(1/n)
\end{eqnarray}
(with $\psi \equiv \Gamma'/\Gamma$ estimated by the Stirling formula, 
$\gamma =$ Euler's constant, and using the polar structure (\ref{DPol}) 
for $Z(\sigma)$).

Prior to using (\ref{EX}) to fix the $R_{-j}$, the argument
covers zeros $\rho$ of a more general (arithmetic) Dirichlet series $L(s)$:
as long as the latter has 
a meromorphic structure and functional equation similar enough to $\zeta(s)$, 
its secondary zeta function $Z(\sigma)$ keeps a double pole at $\sigma=\hf$
\cite{Vl}.
A related but more concrete requirement can be put on the function $N(T)$, 
the number of zeros of $L(s)$ with $0< \Im \rho <T$: 
we ask that for some constants $R_{-2},\ R_{-1}$ and some $\alpha <1$
(all now depending on the chosen $L$-series),
\begin{equation}
\label{NT}
N(T) = 2T \, [ 2 R_{-2} (\log T -1) + R_{-1} ] + \delta N(T), \quad
\delta N(T) = O(T^\alpha) \mbox{ for } T \to +\infty
\end{equation}
(implying $R_{-2} \ge 0$).
If both conditions (\ref{DPol}) and (\ref{NT}) hold, 
then the polar coefficients of $Z(\sigma)$ 
in (\ref{DPol}) have to be precisely the $R_{-j}$ from (\ref{NT}).
All of that is realized in many cases including, but not limited to,
Dedekind zeta functions \cite{La,Ku,Vl} 
and some Dirichlet $L$-functions \cite{Da,Vl,LI2};
see also \cite{Lg} (discussed at end); 
in all those instances, $\delta N(T) = O(\log T)$.
For the corresponding $\lambda_n$, our saddle-point evaluation then always yields: 
either (\ref{NT})~$\Rightarrow$~(\ref{REs}) if all the zeros have $\Re \rho = \hf$ 
-- or the immediately general result (\ref{DA}) otherwise.
\medskip

Like (\ref{DA}) before, (\ref{REs}) can be derived quite rigorously 
but by still another method, previously unknown to us, and written 
for the Riemann zeros by J. Oesterl\'e \cite{O} (private communication). 
We thank him for allowing us to repeat his argument here; 
we actually word it in the more general present setting (and slightly \emph{streamline it}).
When all the zeros lie on the critical line,
first transform the summation (\ref{LDef}) into a Stieltjes integral
(where $\theta(T) = 2 \arctan (1/ \, 2T)$): 
$ \lambda_n = 2 \int_0^\infty [1-\cos n \theta(T)] \,\d N(T) $, 
then integrate by parts:
$ n^{-1} \lambda_n  = 
2 \int_0^\pi \sin n \theta \, N ( \hf \cot {\theta \over 2} ) \,\d \theta $.
Now replace $N(T)$ by its large-$T$ form (\ref{NT}) neglecting $\delta N(T)$ 
and other $O(\theta^{-\alpha} )$ terms: the error is $o(1)$ 
\emph{by the Riemann--Lebesgue lemma}, mainly because \emph{(\ref{NT}) makes
$\delta N(\hf \cot {\theta \over 2})$ integrable over the closed interval} 
$[0,\pi]$; 
hence $ n^{-1} \lambda_n  = \int_0^\pi \sin n \theta \, {1 \over \theta} 
[8 R_{-2} (\log {1 \over \theta} -1) + 4 R_{-1}] \,\d \theta + o(1) $.
Next, change variable: $n \theta = t$, then replace upper $t$-bound 
$n\pi$ by $+\infty$ again with an $o(1)$ error; thus, $ \lambda_n =
{n \int_0^\infty {\sin t \over t} [8 R_{-2} (\log {n \over t} -1) + 4 R_{-1}] 
\,\d t} \pmod{o(n)}$, 
and the last integral evaluates in closed form \cite{GR} to yield (\ref{REs}).
Unfortunately, we do not see how to extend 
this purely real-analytic argument to include the [RH false] case as well.

\subsection*{Recapitulation}

As we ended up with two mutually exclusive large-$n$ behaviors for the $\lambda_n$, 
(\ref{DA}) and (\ref{REs}), together they provide a sharp equivalence result.
For the Riemann zeros, using the explicit values (\ref{EX}):
\medskip

{\bf Theorem} (asymptotic criterion for the Riemann Hypothesis).
\emph{For $n \to +\infty$, the sequence $\lambda_n$ built over the Riemann zeros
follows one of these asymptotic behaviors:}
\begin{eqnarray}
\mbox{\bf - [RH true]} \qquad &\Leftrightarrow& \qquad
\mbox{\emph{tempered}  growth to $+\infty$, as} \nonumber \\
\label{RER}
\lambda _n &\sim& \hf n (\log n - 1+\gamma  - \log 2\pi) \quad \pmod {o(n)} ; \\[6pt]
\mbox{\bf - [RH false]} \qquad &\Leftrightarrow& \qquad
\mbox{\emph{non-tempered} oscillations, as} \nonumber \\
\label{RS}
\lambda _n &\sim& \!\! \sum_{\{ \arg \tau_k >0 \}}
\Bigl( {\tau_k+\mi/2 \over \tau_k-\mi/2} \Bigr) ^n
+ {\rm c.c.} \quad \pmod {o(\e^{\eps n}) \ \forall \eps >0}.
\end{eqnarray}

This comprehensive asymptotic statement \cite{V} is new on the [RH false] side 
to our knowledge, and it also completes some earlier results in the [RH true] case.

Our main end formula ((\ref{RER}), assuming RH) had actually been displayed 
by Keiper \cite[Equation~(37)]{K}, but under a somewhat misleading context: 
in his words, (\ref{RER}) required not just RH but also
``very evenly distributed" zeros,  and was ``much stronger than" RH; 
no details or proofs were ever supplied. All that hardly points toward 
our present conclusion that (\ref{RER}) and RH are strictly equivalent.

Oesterl\'e had a proof of the statement [RH true] $\Rightarrow$ (\ref{RER})
(see above), but he neither published nor even posted his typescript \cite{O}.

In \cite[Cor.~1(c)]{BL}, rather weak exponential lower bounds 
$\lambda_n \ge -c {\rm e}^{\eps n}$ were shown to imply RH; 
the backward assertion [RH true] $\Leftarrow$ (\ref{RER}) is thus also 
\emph{implied} by \cite{BL} (but cannot be inferred therefrom, 
as \cite{BL} never alludes to asymptotics regarding the $\lambda_n$).

Our saddle-point approach also handles 
both cases (\ref{RER})--(\ref{RS}) at once for the first time.
\medskip

Numerical data \cite{K,M1} agree well with (\ref{RER}) for $n<7000$ 
(and even better in the mean if we add the contribution
like (\ref{REs}) but from the next pole of $I(\sigma)$,
\begin{equation}
\delta \lambda _n = (-1)^n 2n \Res_{\sigma=0} I(\sigma) = 2Z(0) =  +7/4
\end{equation}
\cite[Equation~(41)]{Vz}, although this correction should not count asymptotically, 
dominated as it seems by oscillatory terms). 
Yet the above numerical agreement is inconclusive regarding the Riemann Hypothesis: 
any currently possible violation of RH would yield a deviation from (\ref{RER}) 
detectable only at \emph{much higher} $n$ (see end of Sect.~4).

\section{An even more sensitive sequence.}

A slightly stronger difference of behavior follows for
the special linear combinations (\ref{SN}) below 
of the coefficients $\eta_n$ themselves (defined by (\ref{SC}) above).
 
Indeed, the definition 
$\Xi(s)=  \Gamma (s/2) \pi^{-s/2}  s(s-1) \zeta(s)$
substituted into (\ref{LId}) readily yields a decomposition 
$\lambda_n=S_n+\overline S_n \,$, where \cite[thm~2]{BL} 
\begin{equation}
\label{SN}
S_n = - \sum_{j=1}^n {n \choose j} \, \eta_{j-1} 
\end{equation} 
is the contribution of $(s-1)\zeta(s)$, and
\begin{equation}
\label{TN}
\overline S_n = 1 - {\log 4\pi + \gamma \over 2} \, n + \hat S_n \, ,\quad 
\mbox{with} \quad 
\hat S_n = \sum _{j=2}^n {n \choose j} (-1)^j(1-2^{-j}) \, \zeta(j) ,
\end{equation}
is the contribution of the remaining (more explicit) factor.
The large-$n$ behavior of the sum $\hat S_n$ is also computable, 
by the same route we followed from (\ref{LZS}) to (\ref{REs}) through (\ref{LZI}). 
First,
\begin{equation}
\hat S_n =
{(-1)^n n! \over 2\pi\mi} \oint_{C'} J(\sigma) \,\d \sigma, \qquad J(\sigma) = 
{\Gamma(\sigma -n) \over \Gamma(\sigma +1)} (1-2^{-\sigma}) \, \zeta(\sigma) ,
\end{equation}
integrated around the poles $\sigma=2,\ldots,n$ of $J(\sigma)$;
hence for $n \to +\infty$, $\hat S_n$ is asymptotically
\begin{equation}
(-1)^{n-1} n! \Res_{\sigma=1} J(\sigma) = \hf n [ \psi(n) + \log 2 - 1 + 2 \gamma ] ,
\end{equation}
now (mod $o(n^{-N})\ \forall N >0$) because $J(\sigma)$ has no further singularities;
so that finally, 
using the Stirling expansion for $\psi(n)$ in terms of the Bernoulli numbers $B_{2k}$,
\begin{equation}
\label{TS}
\overline S_n \sim \hf n (\log n - 1+\gamma  - \log 2\pi) + {3 \over 4} 
- \sum_{k=1}^\infty {B_{2k} \over 4k} n^{1-2k} \quad (n \to \infty) \quad
\mbox{\emph{unconditionally}} .
\end{equation}

This at once confirms two empirical conjectures made by Ma\'slanka 
\cite[Equations~(2.5), (2.8)]{M1}:
mod~$o(n)$, the sequence $\{ \overline S_n \}$ expresses the ``trend" (\ref{RER}) 
obeyed by the sequence $\{\lambda_n\}$ \emph{under [RH true]};
and to all orders in $n$, 
$\{ \overline S_n \}$ has the asymptotic expansion (\ref{TS}).
But whether RH holds or not it makes sense to withdraw 
the fixed $\{\overline S_n\}$-contribution from the previous formulae 
(\ref{RER}) and (\ref{RS}), as we did in \cite{V} to find:
\begin{equation}
\label{SNT}
S_n = o(n) \qquad \qquad \qquad \qquad \qquad \qquad \qquad \qquad \qquad 
\qquad \qquad \qquad \qquad [{\rm RH\ true}]
\end{equation}
(a case further discussed in \cite{M2,S,C2}), versus
\begin{equation}
\label{SLN}
S_n \sim \lambda_n \sim  \sum_{\{ \arg \tau_k >0 \}}
\Bigl( {\tau_k+\mi/2 \over \tau_k-\mi/2} \Bigr) ^n
{}+ {\rm c.c.} \quad \pmod {o(\e^{\eps n}) \ \forall \eps >0}
\qquad \quad [{\rm RH\ false}],
\end{equation}
which gives oscillations that blow up \emph{exponentially} with $n$.

Still, in absolute size, any contribution from (\ref{SLN}) will stay 
considerably smaller than (\ref{SNT}) (the background from the set of real $\tau_j$)
up to $n \approx \min_{\{ \arg \tau_k >0 \}} \{ |\Im 1/\tau_k|^{-1} \}$:
i.e., $S_n$ can only reliably signal a zero violating RH up to a height 
$| \Im \rho | \lesssim \sqrt{n/2}$; 
so vice-versa, since such zeros are now known to require
$| \Im \rho | \gtrsim 10^9$, they could only be detected by $S_n$ 
for $n \gtrsim 10^{18}$ (see also \cite{BPY,O}, \cite[p.~5]{Lg}). 

Our asymptotic criteria may thus not surpass others in practical sensitivity,
but their sharpness is theoretically interesting.
For instance, Li's criterion ${\lambda_n >0} \ (\forall n)$ is now strengthened
in that, beyond a finite $n$-range, it refers to ever larger amplitudes either way, 
i.e., delicate borderline situations such as $\lambda_n \to 0$ are ruled out.

\section{More general cases: a summary}

For the more general zeta-type functions satisfying (\ref{DPol}) and (\ref{NT}),
the asymptotic alternative for the associated $\lambda_n$ is:
either (\ref{REs}) if the generalized Riemann Hypothesis (GRH) holds, 
or (\ref{RS}) otherwise. 
In the former case, the connection (\ref{NT}) $\Rightarrow$ (\ref{REs}) often 
makes the asymptotic form of $\lambda_n$ fully explicit without any calculation,
and it also ensures $R_{-2} \ge 0$
in line with the corresponding generalized Li's criterion \cite{BL,Lg,LI2,LI3}.

The $\lambda_n$ have also been generalized to $L$-functions defined by 
Hecke operators for the congruence subgroup $\Gamma_0(N)$
\cite{LI2}\cite[specially remark~5.4]{LI3}.

More recently and in a broader setting 
(the $\lambda_n$ for automorphic $L$-functions), Lagarias presented 
an alternative approach to estimate the $\lambda_n$ with greater accuracy, 
mod $O(\sqrt n \log n)$ under GRH \cite{Lg}: 
in the notations of (\ref{SN})--(\ref{TN}), he directly proves that 
$\overline S_n$ obeys (\ref{REs}) mod $O(1)$ unconditionally (thm~5.1), 
then that $S_n = O(\sqrt n \log n)$ under GRH (thm~6.1).
We note that his leading $\lambda_n$-behavior remains tied to 
the large-$T$ behavior of a counting function according to the rule
(\ref{NT}) $\Rightarrow$ (\ref{REs}) \cite[Equations~(2.11--12),(1.18),(5.2)]{Lg}.


\begin{thebibliography}{99}

\bibitem{BPY} Biane, P., Pitman, J. and Yor, M.: Probability laws related to 
the Jacobi theta and Riemann zeta functions, 
{\it Bull. Amer. Math. Soc.\/} {\bf 38} (2001), 435--465 [Sect.~2.3].

\bibitem{BL} Bombieri, E. and Lagarias, J.C.:
Complements to Li's criterion for the Riemann Hypothesis,
{\it J. Number Theory\/} {\bf 77} (1999), 274--287.

\bibitem{C1} Coffey, M.W.: 
Relations and positivity results for the derivatives of 
the Riemann $\xi$ function,
{\it J. Comput. Appl. Math.\/} {\bf 166} (2004), 525--534.

\bibitem{C2} Coffey, M.W.: New results concerning power series expansions of
the Riemann xi function and the Li/Keiper constants, preprint (Jan. 2005);
Toward verification of the Riemann Hypothesis:
application of the Li criterion, {\it Math. Phys. Anal. Geom.\/} 
{\bf 8} (2005) 211--255.

\bibitem{Da} Davenport, H.: {\it Multiplicative Number Theory\/},
3rd ed., revised by H.L. Montgomery, Graduate Texts in Mathematics {\bf 74},
Springer-Verlag (2000) [chap.~16].

\bibitem{Di} Dingle, R.B.: 
{\it Asymptotic Expansions: Their Derivation and Interpretation\/}, 
Academic Press (1973).

\bibitem{Ed} Edwards, H.M.: {\it Riemann's Zeta Function\/}, 
Academic Press (1974) [Sect.~1.10].

\bibitem{Er} Erd\'elyi, A.: {\it  Asymptotic Expansions\/}, Dover (1956)
[Sect.~2.5].

\bibitem{GR} Gradshteyn, I.S. and Ryzhik, I.M.: 
{\it Table of Integrals, series and products\/}, 5th ed., A. Jeffrey ed.,
Academic Press (1994) [Equations (3.721(1)) p.~444 and (4.421(1)) p.~626].

\bibitem{H} Hashimoto, Y.: Euler constants of Euler products,
{\it J. Ramanujan Math. Soc.\/} {\bf 19} (2004), 1--14.

\bibitem{I} Israilov, M.I.: 
On the Laurent expansion of the Riemann zeta-function,
{\it Proc. Steklov Inst. Math.\/} {\bf 4} (1983), 105--112 
[Russian: {\it Trudy Mat. Inst. i. Steklova\/} {\bf 158} (1981), 98--104].

\bibitem{K} Keiper, J.B.:
Power series expansions of Riemann's $\xi$ function,
{\it Math. Comput.\/} {\bf 58} (1992), 765--773.

\bibitem{Ku} Kurokawa, N.: Parabolic components of zeta functions,
{\it Proc. Japan Acad.\/} {\bf 64}, Ser.~A (1988), 21--24;
Special values of Selberg zeta functions, in:
{\it Algebraic K-theory and algebraic number theory\/}
(Proceedings, Honolulu 1987), M.R. Stein and R. Keith Dennis eds.,
{\it Contemp. Math.\/} {\bf 83}, Amer. Math. Soc. (1989), pp.~133--149.

\bibitem{Lg} Lagarias, J.C.: 
Li coefficients for automorphic $L$-functions, 
{\it Ann. Inst. Fourier, Grenoble\/} (2006, to appear) {\tt [math.NT/0404394 v4]}.

\bibitem{La} Landau, E.: {\it Einf\"uhrung in die elementare und analytische 
Theorie der algebraischen Zahlen und der Ideale\/}, Chelsea, New York (1949)
[Satz~173 p.~89].

\bibitem{Le} Lehmer, D.H.: 
The sum of like powers of the zeros of the Riemann zeta function, 
{\it Math. Comput.\/} {\bf 50} (1988), 265--273.

\bibitem{LI1} Li, X.-J.:
The positivity of a sequence of numbers and the Riemann Hypothesis,
{\it J. Number Theory\/} {\bf 65} (1997), 325--333.

\bibitem{LI2} Li, X.-J.: 
Explicit formulas for Dirichlet and Hecke $L$-functions,
{\it Illinois J. Math.\/} {\bf 48} (2004), 491--503.

\bibitem{LI3} Li, X.-J.: An arithmetic formula for certain coefficients of
the Euler product of Hecke polynomials, 
{\it J. Number Theory\/} {\bf 113} (2005), 175--200.

\bibitem{M1} Ma\'slanka, K.: 
Effective method of computing Li's coefficients and their properties, 
{\it Experiment. Math.\/} (to appear) {\tt [math.NT/0402168 v5]}.

\bibitem{M2} Ma\'slanka, K.:
An explicit formula relating Stieltjes constants and Li's numbers, 
preprint {\tt [math.NT/0406312 v2]}.

\bibitem{M} Matsuoka, Y.: 
A note on the relation between generalized Euler constants 
and the zeros of the Riemann zeta function,
{\it J. Fac. Educ. Shinshu Univ.\/} {\bf 53} (1985), 81--82;
A~sequence associated with the zeros of the Riemann zeta function, 
{\it Tsukuba J. Math.\/} {\bf 10} (1986), 249--254.

\bibitem{O} Oesterl\'e, J.: 
R\'egions sans z\'eros de la fonction z\^eta de Riemann,
typescript (2000, revised 2001, uncirculated).

\bibitem{S} Smith, W.D.: 
A ``good" problem equivalent to the Riemann hypothesis,
e-print on http://www.math.temple.edu/$\scriptstyle\sim$ wds/homepage/works.html
(2005 version, unpublished).

\bibitem{VP} de la Vall\'ee Poussin, C.-J.: 
Recherches analytiques sur la th\'eorie des nombres premiers~I,
{\it Ann. Soc. Sci. Bruxelles\/} {\bf 20} (1896), 183--256 [p.~251].

\bibitem{Vz} Voros, A.: Zeta functions for the Riemann zeros,
{\it Ann. Inst. Fourier, Grenoble\/} {\bf 53} (2003), 665--699;
erratum: {\bf 54} (2004), 1139.

\bibitem{Vl} Voros, A.:
Zeta functions over zeros of general zeta and $L$-functions, in:
{\it Zeta Functions, Topology and Quantum Physics\/} (Proceedings, Osaka,
March 2003), T.~Aoki, S.~Kanemitsu, M.~Nakahara and Y.~Ohno eds., 
Developments in Mathematics {\bf 14}, Springer-Verlag (2005), pp.~171--196.

\bibitem{V} Voros, A.: 
A sharpening of Li's criterion for the Riemann Hypothesis,
preprint (Saclay-T04/040 April 2004, unpublished) {\tt [math.NT/0404213 v2]}.

\end{thebibliography}
\end{document}